\documentclass[11pt,reqno]{amsart}
\calclayout
\usepackage{amssymb}
\newif\ifaddpics\addpicstrue
\ifaddpics\usepackage{graphicx}\fi
\makeatletter
\def\swappedhead@plain#1#2#3{%
  \thmnumber{\@upn{\mdseries #2}}\thmname{\@ifnotempty{#2}{. }#1}%
  \thmnote{ \textmd{\upshape(#3)}}}
\makeatother
\theoremstyle{plain}

\swapnumbers
\newtheorem{thm}[subsection]{Theorem}
\newtheorem{lem}[subsection]{Lemma}
\newtheorem{prop}[subsection]{Proposition}

\numberwithin{equation}{section}

\newcommand{\thismonth}{\ifcase\month\or
  January\or February\or March\or April\or May\or June\or
  July\or August\or September\or October\or November\or December\fi
  \space\number\year}
\makeatletter
\newcommand{\low}{\@ifnextchar^{}{^{\vphantom x}}}
\newcommand{\high}{\@ifnextchar_{}{_{\vphantom I}}}
\makeatother
\DeclareSymbolFont{script}{U}{eus}{m}{n}
\DeclareSymbolFontAlphabet{\mathscr}{script}
\DeclareMathSymbol{\EuWedge}{0}{script}{"5E}
\DeclareMathAlphabet{\mathrmsl}{OT1}{cmr}{m}{sl}
\newcommand{\bbsymb}[2]{\newcommand{#1}{{\mathbb{#2}}}}

\newcommand{\oper}[3][n]{\newcommand{#2}{\mathop
  {\mathrm{#3}\null}\ifx n#1\nolimits\else\limits\fi}}
\newcommand{\rsoper}[3][n]{\newcommand{#2}{\mathop
  {\mathrmsl{#3}\null}\ifx n#1\nolimits\else\limits\fi}}
\bbsymb\C{C} \bbsymb\F{F} \bbsymb\IH{H}\bbsymb\N{N} \bbsymb\IP{P}
\bbsymb\Q{Q} \bbsymb\R{R} \bbsymb\U{U} \bbsymb\V{V} \bbsymb\W{W}
\bbsymb\Z{Z}
 \bbsymb\E{E}

\renewcommand{\phi}{\varphi}

\newcommand{\ve}{\varepsilon}

\newcommand{\hf}{{\mbox{\small$\frac{1}{2}$}}}

\newcommand{\oZ}{\overline{Z}}
\newcommand{\oD}{\overline{D}}
\newcommand{\tx}{\tilde{x}}
\newcommand{\ty}{\tilde{y}}

\newcommand{\cH}{\mathcal H}
\newcommand{\oH}{\overline{\mathcal H}}
\newcommand{\rd}{\mbox{d}}

\newcommand{\sign}{\mbox{sgn}}
\renewcommand{\geq}{\geqslant} \renewcommand{\leq}{\leqslant}
\theoremstyle{definition}

\newcommand{\del}{\partial}                 

\begin{document}
\title[Continued fractions and Einstein manifolds]{Continued fractions
and Einstein manifolds\\ of infinite topological type}
\author{David M. J. Calderbank}
\address{School of Mathematics\\
University of Edinburgh\\ King's Buildings, Mayfield Road\\ Edinburgh
EH9 3JZ\\ Scotland.}
\email{dc511@york.ac.uk}
\curraddr{Department of Mathematics\\ University of York\\ Heslington\\
York YO10 5DD\\ England.}
\author{Michael A. Singer}
\email{m.singer@ed.ac.uk}

\begin{abstract}  We present a construction of complete self-dual Einstein
metrics of negative scalar curvature on an uncountable family of manifolds of
infinite topological type.
\end{abstract}

\maketitle
\section{Introduction}

\subsection{Summary}

Let $\alpha$ be an irrational number, $0<\alpha<1$, and consider the {\em 
modified continued fraction expansion} of $\alpha$,
\begin{equation}\label{e1.9.6.5}
\alpha = \cfrac{1}{e_1-\cfrac{1}{e_2- \cfrac{1}{e_3-\cdots}}},
\end{equation}
with $e_j\geq 2$ for all $j$.  We associate to $\alpha$ a noncompact
$4$-manifold $M_\alpha$, which is connected and simply connected, but has
$b_2(M_\alpha)= \infty$, $H_2(M_\alpha)$ being generated by an infinite
sequence of embedded $2$-spheres $S_j$, with
\begin{equation}\label{e2.9.6.5}
S_j \cdot S_{j+1} = -1,\;  S_j\cdot S_j = e_j,\;
S_j\cdot S_k = 0\mbox{ for }|j-k|>1.
\end{equation}
Then, subject to the uniform bound
$3\leq e_j\leq N$, we construct a complete self-dual Einstein metric
$g_\alpha$, of negative scalar curvature, on $M_\alpha$. The key point about
the construction is that $M_\alpha$ and $g_\alpha$ are toric: there is a smooth
action of $T^2$ on $M_\alpha$ which preserves $g_\alpha$.

The work in this paper complements that in our previous papers \cite{CaSi:emcs}
and \cite{CaSi:temco}.  Indeed, in \cite{CaSi:emcs}, we associated to each {\em
rational} number $\alpha$, $0<\alpha<1$, a noncompact $4$-manifold $M_\alpha$
and, subject to the condition $e_j\geq 3$, we constructed a complete SDE metric
$g_\alpha$ on $M_\alpha$, where, as above, $M_\alpha$ is connected and simply
connected, but now $H_2(M_\alpha)$ is generated by a \emph{finite} sequence of
spheres $S_j$ which satisfy \eqref{e2.9.6.5}.

We note a parallel in hyperk\"ahler geometry: the $A_n$-gravitational
instantons constructed by the Gibbons--Hawking Ansatz are analogous to the case
that $\alpha$ is rational, while the hyperk\"ahler manifolds of infinite
topological type in \cite{AKL:rfitt} correspond to $\alpha$ being irrational.
But there are very important differences: in the hyperk\"ahler case, the
self-intersections $e_j$ all have to be $2$, which is complementary to the
condition $e_j\geq 3$ that we impose here. On the other hand, despite the
uniform bound $e_j \leq N$, we obtain uncountably many non-diffeomorphic SDE
manifolds in this way---for this, it is enough to allow the $e_j$ to
take only the values $3$ and $4$, for example.

The plan of this paper is as follows.  In \S\ref{s2.14.6.5} we gather some
elementary facts about the continued fraction expansion \eqref{e1.9.6.5}. The
only result that may be new here is Theorem~\ref{t1.3.5.5}; we are grateful to
Chris Smyth for assisting us with its proof.  In \S\ref{s10.14.6.5}, we give
the construction of $M_\alpha$. This is a straightforward extension of the work
in \cite{CaSi:emcs}, which, as we have indicated, corresponds to the case that
$\alpha$ is rational. (That work, in turn, rests on the combinatorial
description of toric $4$-manifolds due to Orlik and Raymond \cite{OrRa:at4}.)
In \S\ref{s11.14.6.5} we write down a SDE metric $g_\alpha$ on $M_\alpha$.
This is defined initially on a dense open subset $U\subset M_\alpha$ (the set
on which the $T^2$-action is free) but at the end of the section we show that
$g_\alpha$ extends smoothly to the whole manifold. The argument is given in
detail partly to make this paper more self-contained, and also to clarify one
point omitted from \cite{CaSi:emcs}\footnote{For readers of that paper, we did
not show in \cite[\S 5]{CaSi:emcs} that $\lim_{\rho\to 0}
\rho^{-1}|\det\Phi(\rho,\eta)| > 0$, and this is needed for the smooth
extension of the metric.}. In \S\ref{s1.5.7.5} we give some technical estimates
on the functions which enter the definition of $g_\alpha$. These are needed for
the smooth extension of $g_\alpha$ from $U$ to $M_\alpha$, and also pave the
way for the proof that $g_\alpha$ is complete in \S\ref{s1.9.8.5}.

For the reader familiar with \cite{CaSi:emcs}, we make some remarks about the
extension of the construction of $g_\alpha$ from finite to infinite continued
fractions. From the point of view of \cite[\S4]{CaSi:emcs}, one would like to
allow $k$ to go to $\infty$ in
\begin{equation*}
\sqrt{\rho}F(\rho,\eta) = \sum_{j=0}^{k+1}w_j \sqrt{\rho^2 +(\eta- y_j)^2},
\end{equation*}
where
\begin{equation*}
w_j = m_{j+1}-m_j, \quad y_j = \frac{n_{j+1}-n_j}{m_{j+1}-m_j}
\end{equation*}
and the $n_j/m_j$ are the continued fraction approximants to
$\alpha$. But the sequence $w_j$ increases rapidly with $j$, while the
$y_j$ also converge to $\alpha$, so the status of such a limit is
unclear.  However, we noted in \cite[\S5]{CaSi:emcs} that the
above sum represents an eigenfunction of the hyperbolic laplacian with
boundary data equal to $\eta\mapsto m_j \eta -n_j$ for $y_j \leq \eta
\leq y_{j-1}$ and this makes perfectly good sense also for infinite
continued fractions.  This observation is really the key to the
construction in this paper.

\subsection{Notation} \label{s1.14.6.5} Denote by $\cH^2$ the
hyperbolic plane, by $\oH{}^2$ its conformal compactification. We shall
always identify $\cH^2$ with the upper half-plane $\{(x,y)\in \R^2:
y>0\}$; then the hyperbolic metric is given by $(\rd x^2 +\rd
y^2)/y^2$ and $\oH{}^2 =\{(x,y): y \geq 0\}\cup\{\infty\}$.  Note that
half-space coordinates near $\infty$ can be defined by
$$
(\tx, \ty) = \left(-\frac{x}{x^2+y^2}, \frac{y}{x^2+y^2}\right).
$$
Further, denote by $T^2$ the standard $2$-torus, identified with
$\R^2/(2\pi \Z)^2$.  We shall write $z = (z_1,z_2)$ for standard
linear coordinates on $\R^2$. The circle-subgroup generated by
$m\del_{z_1} + n\del_{z_2}$ (for $m,n\in\Z$) will be denoted by $S^1_{(m,n)}$.
We shall denote by $\ve$ the standard skew form
\begin{equation}\label{e1.9.8.5}
\ve(z',z'') = \det(z',z'') = z_1'z''_2 - z'_2z''_1.
\end{equation}

\subsection{Acknowledgement} We thank Chris Smyth for useful
conversations on continued fractions and for the proof of
Theorem~\ref{t1.3.5.5}.  We also thank Jim Wright for useful conversations.

\section{Continued fractions and toric $4$-manifolds}
\label{s2.14.6.5}
\subsection{Continued fractions}
\label{s3.14.6.5}
Let us begin with our irrational number $\alpha$ and its continued fraction
expansion \eqref{e1.9.6.5}.  Set
\begin{equation}\label{e2.2.5.5}
(m_0,n_0) = (0,-1),\quad (m_0,n_0) = (1,0),
\end{equation}
and
\begin{equation}\label{e3.2.5.5}
(m_j,n_j)\mbox{ coprime with }
\frac{n_j}{m_j} = 
 \cfrac{1}{e_1-\cfrac{1}{e_2- \cdots \cfrac{1}{e_{j-1}}}}.
\end{equation}
The most important properties of this sequence of pairs $(m_j,n_j)$
are summarized as follows.
\begin{lem} \label{l2.14.6.5}
For each $j\geq 1$,
\begin{equation}\label{e6.2.5.5}
(m_{j+1},n_{j+1}) = e_j(m_j,n_j) - (m_{j-1},n_{j-1}),
\end{equation}
and for $j\geq 0$,
\begin{equation}\label{e7.2.5.5}
m_{j}n_{j+1} - m_{j+1}n_{j} = 1.
\end{equation}
\end{lem}
\begin{proof}
If 
\begin{equation}\label{e5.2.5.5}
M_j = \begin{pmatrix} 0& 1 \cr -1 & e_j\end{pmatrix}
\end{equation}
then one has, for each $j$:
\begin{align}
\begin{pmatrix} n_{j-1} \cr m_{j-1} \end{pmatrix} &=
M_1M_2\cdots M_{j-1}
\begin{pmatrix} -1\cr 0 \end{pmatrix}; \label{e1.14.6.5} \\
\begin{pmatrix} n_j \cr m_j \end{pmatrix} &=
M_1M_2\cdots M_{j-1}
\begin{pmatrix} 0\cr 1 \end{pmatrix};\label{e4.2.5.5} \\
\begin{pmatrix} n_{j+1} \cr m_{j+1} \end{pmatrix} &=
M_1M_2\cdots M_{j-1}
\begin{pmatrix} 1\cr e_j \end{pmatrix}.  \label{e2.5.6.5} 
\end{align}
These are all essentially equivalent to each other and are easily
proved by induction. Combining these, we have
\begin{equation}\label{e1.5.6.5}
\begin{pmatrix} n_{j+1} & n_{j} \cr
 m_{j+1} & m_{j} \end{pmatrix} = M_1M_2\cdots M_{j-1}
 \begin{pmatrix} 1 & 0 \cr e_j & 1 \end{pmatrix},\quad
\begin{pmatrix} n_{j+1} & n_{j-1} \cr
 m_{j+1} & m_{j-1} \end{pmatrix} = M_1M_2\cdots M_{j-1}
 \begin{pmatrix} 1 & -1 \cr e_j & 0 \end{pmatrix},
 \end{equation}
and, taking determinants, we obtain \eqref{e7.2.5.5} and the identity
\begin{equation}\label{e8.2.5.5}
m_{j-1}n_{j+1} - m_{j+1}n_{j-1} = e_j.
\end{equation}
From \eqref{e7.2.5.5} (with $j$ replaced by $j-1$) it follows that
$(m_{j-1},n_{j-1})$ and $(m_j,n_j)$ form a $\Z$-basis of
$\Z\oplus\Z$, so that $(m_{j+1},n_{j+1})$ is an integer linear combination
of these vectors. The coefficients are
determined as in \eqref{e6.2.5.5} by using the identities
\eqref{e7.2.5.5} and \eqref{e8.2.5.5}.
\end{proof}

\subsection{Definition} Set 
\begin{equation}\label{e9.2.5.5}
a_j = \frac{n_{j+1}-n_j}{m_{j+1}-m_j},\quad b_j = \frac{1}{m_{j+1}-m_j}.
\end{equation}
These are well-defined by the following lemma.
\begin{lem} \label{l1.14.6.5}
\textup{(i)} For all $j\geq 1$, $m_{j+1}> m_j$, $n_{j+1}>n_j$, and $(n_j/m_j)$
is strictly increasing with limit $\alpha$. Hence the sequences $(a_j)$ and
$(b_j)$ are positive. Furthermore, $\lim_{j\to\infty} a_j=\alpha$.

\noindent\textup{(ii)} Suppose $e_j\geq 3$ for all $j$.  Then $m_{j+1} > \phi^2
m_j$ and $n_{j+1} > \phi^2 n_j$, where $\phi=(1+\sqrt 5)/2>1$ is the golden
ratio.  Furthermore, $(a_j)$ and $(b_j)$ are strictly decreasing, and
$\lim_{j\to\infty} b_j=0$.
\end{lem}
\begin{proof} (i) To prove that the sequences $(m_j)$ and $(n_j)$ are
strictly increasing, use \eqref{e6.2.5.5} and induction on $j$. In particular
$m_j>0$ for $j>0$. The monotonicity of the sequence $n_j/m_j$ now follows from
the identity
\begin{equation}\label{e2.14.6.5}
\frac{n_{j+1}}{m_{j+1}} - \frac{n_{j}}{m_{j}}  
= \frac{1}{m_jm_{j+1}},
\end{equation}
where we have used \eqref{e7.2.5.5}. The fact that $n_j/m_j\to\alpha$ as
$j\to\infty$ is standard. Now note that
\begin{equation}\label{e3.3.5.5}
a_j -\frac{n_j}{m_j} = \frac{1}{m_j(m_{j+1}-m_j)}>0
\end{equation}
so that $\lim_{j\to\infty} a_j = \lim_{j\to\infty} n_j/m_j = \alpha$ as
required.

(ii) We again use~\eqref{e6.2.5.5} and induction on $j$: e.g., if $m_j>\phi^2
m_{j-1}$, then
\begin{equation*}
m_{j+1} = e_j m_j - m_{j-1} > (3-\phi^{-2})m_j =\phi^2 m_j,
\end{equation*}
since $\phi^2=(3+\sqrt 5)/2$.  Now from \eqref{e6.2.5.5} and the nonnegativity
of $(m_j)$,
\begin{equation*}
m_{j+1} -m_j = (e_j-1)m_j - m_{j-1} \geq (e_j-1)(m_j -m_{j-1}).
\end{equation*}
Since $e_j-1>1$, the statements about $(b_j)$ are immediate. On the other hand,
\begin{equation}\label{e1.3.5.5}
a_{j-1} - a_{j} = \frac{e_j-2}{(m_{j+1}-m_j)(m_j-m_{j-1})}>0
\end{equation}
which gives the monotonicity of the $a_j$.
\end{proof}

\subsection{Assumption}
\label{s5.14.6.5} From now on, suppose that for all $j$,  $3\leq e_j
\leq N$ for some $N>0$. 

\subsection{Definition} The {\em envelope} $\eta_\alpha$ of $\alpha$
is defined as
\begin{equation}\label{e4.3.5.5}
  \eta_{\alpha}(x) = \begin{cases} m_jx -n_j & 
   \mbox{for }x\in [a_j,a_{j-1}],\\
      0& \mbox{for }x\leq \alpha.
\end{cases}
\end{equation}
(Here we set $a_{-1}= \infty$, so that $\eta_{\alpha}(x) =1$ for $x\geq a_0 =
1$.)
\begin{prop}\label{eta-prop} $\eta_\alpha$ is continuous, and for $x>\alpha$ it
is strictly increasing \textup(hence positive\textup), concave \textup(i.e.,
$\eta''_\alpha \leq 0$ in the sense of distributions\textup) and is the linear
interpolant of the points $(a_j,b_j)$ \textup(so that $\eta_\alpha(a_j)=b_j$
for all $j$\textup).
\end{prop}
\begin{proof} Trivial, given the previous results.\end{proof}

\subsection{Example}
\label{s6.14.6.5} Suppose $e_j=3$ for all $j$. Then
\begin{equation*}
m_j = n_{j+1} =  \frac{\phi^{2j} - \phi^{-{2j}}}{\sqrt{5}},\quad
a_j=\frac{\phi^{2j-1}+\phi^{-2j+1}}{\phi^{2j+1}+\phi^{-2j-1}}, \quad
b_j=\frac{\sqrt 5}{\phi^{2j+1}+\phi^{-2j-1}},
\end{equation*}
where $\phi = (1 +\sqrt{5})/2$ as above. Then $\alpha
=\phi^{-2}=0.381966\ldots$ and
\begin{equation*} 
\frac{a_j - \phi^{-2}}{b_j^2} = \frac{1}{\sqrt{5}}
\biggl(\frac{\phi^{4j+2}+2+\phi^{-4j-2}}{\phi^{4j+2}+1}\biggr).
\end{equation*}
In particular, by Proposition~\ref{eta-prop},
\begin{equation} \label{e2.9.8.5}
\eta_{\alpha}(x) \lesssim \sqrt{\sqrt{5}(x-\phi^{-2})}
\end{equation}
for $x-\phi^{-2}$ small.  The two functions in \eqref{e2.9.8.5} are shown in
Figure~\ref{fig1}.
\ifaddpics
\begin{figure}[ht]
\begin{center}
\includegraphics[width=.7\textwidth]{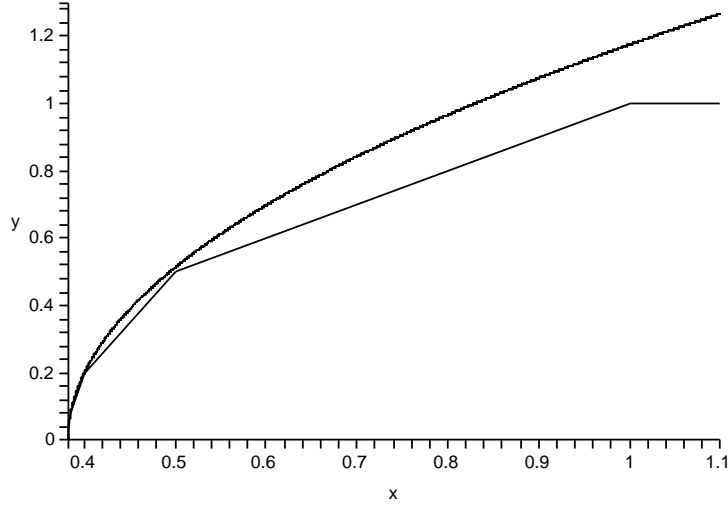}
\caption{The envelope and the square-root function with $\alpha=\phi^{-2}$}
\label{fig1}
\end{center}
\end{figure}
\fi
\smallbreak

We shall now show that the behaviour of $\eta_\alpha(x)$ is bounded by a
multiple of $\sqrt{x-\alpha}$ in general---we are indebted to Chris Smyth for
the proof of this fact.

\begin{thm} \label{t1.3.5.5}
There is a  constant $\Omega>0$ such that
\begin{equation}\label{eta-bounds}
\eta_\alpha(x) \leq  \Omega\sqrt{x-\alpha} \mbox{ for }x\geq\alpha.
\end{equation}
\end{thm}
\begin{proof} This is a sequence of elementary deductions, based on
Lemmas~\ref{l2.14.6.5} and \ref{l1.14.6.5}.
Start by substituting 
$m_{j+1} =e_jm_j -m_{j-1}$ into  
\eqref{e1.3.5.5}, to get
\begin{equation}\label{11.19.12}
a_{j-1} - a_{j} = \left[\frac{e_{j}-2}{(1 - m_{j-1}/m_{j})
(e_{j} -1 - m_{j-1}/m_{j})}\right]\frac{1}{m_{j}^2}.
\end{equation}
Because $0 < m_{j-1}/m_{j} <\phi^{-2}<1/2$ and $e_j\geq 3$ it is elementary to
bound the quantity in square brackets between $1/2$ and $2$. Hence
\begin{equation}\label{12.19.12}
\frac{1}{2m_{j}^2}< a_{j-1} - a_{j}  < \frac{2}{m_{j}^2}.
\end{equation}
Now
\begin{equation}\label{e4.3.6.5}
a_n- \alpha = \sum_{j=n+1}^\infty (a_{j-1} - a_{j})
\end{equation}
and since $m_j > 2m_{j-1}$ for all $j$, we obtain
\begin{equation}\label{13.19.12}
\frac{1}{2m_{n+1}^2} < a_n -\alpha
< \frac{2}{m_{n+1}^2}\sum_{j=0}^\infty \phi^{-4j} <\frac{3}{m_{n+1}^2}.
\end{equation}
(The lower bound is the trivial one coming from the first term in the sum.)

Since $b_{j-1} - b_{j} = m_j(a_{j-1} - a_{j})$ we find from \eqref{12.19.12}
\begin{equation}\label{14.19.12}
\frac{1}{2m_{j+1}}< b_{j-1} - b_{j} < \frac{2}{m_{j+1}};
\end{equation}
summing as before, we obtain
\begin{equation}\label{15.19.12}
\frac{1}{2m_{n+1}} < b_n < \frac{2}{m_{n+1}}\sum_{j=0}^\infty \phi^{-2j}
< \frac{4}{m_{n+1}}.
\end{equation}
Combining \eqref{13.19.12} and \eqref{15.19.12}  in the obvious way we 
obtain
\begin{equation}\label{16.19.12}
\frac{1}{2\sqrt 3}\sqrt{a_n-\alpha}  <  b_n  <  4\sqrt 2\sqrt{a_n - \alpha}.
\end{equation}
Since $\eta_\alpha$ is piecewise linear, $\eta_\alpha(a_n) =b_n$, and the
square-root function is concave, the bound~\eqref{eta-bounds} follows at once.
\end{proof}

The following technical result will be needed in \S\ref{s1.9.8.5}. In order to
state it, set
\begin{equation}\label{e1.10.8.5}
\mu(x) = m_j\mbox{ if $a_j < x < a_{j-1}$},\quad
\nu(x) = n_j\mbox{ if $a_j < x < a_{j-1}$},
\end{equation}
and
\begin{equation}\label{e2.10.8.5}
D(x_1,x_2) = (\mu(x_1)\nu(x_2) -\mu(x_2)\nu(x_1))(x_1-x_2).
\end{equation}

\begin{lem}\label{l1.11.8.5}
If $3 \leq e_j \leq N$ for all $j$, then
\begin{equation}\label{e3.10.8.5}
D(\alpha+ x, \alpha+ x\eta) \geq  x(1-\eta)\mbox{ for }x\in (0,1-\alpha)),
\eta \in(0,(4N)^{-2}). 
\end{equation}
\end{lem}
\begin{proof}
We have
\begin{equation} \label{e11.5.7.5}
D(\alpha+x,\alpha+ x\eta) = [\mu(\alpha + x)\nu(\alpha + x\eta) -
\mu(\alpha + x\eta)\nu(\alpha + x)]x(1-\eta);
\end{equation}
since the quantity in square brackets is $\geq 1$ if $\alpha+x$ and
$\alpha +x\eta$ are in disjoint intervals $[a_j,a_{j-1}]$,
it is enough to show that if $x$ and $\eta$ are as in
\eqref{e3.10.8.5} and 
\begin{equation}\label{e1.11.8.5}
a_j<\alpha+ x< a_{j-1},
\end{equation} 
then $\alpha + x\eta < a_j$. Now if \eqref{e1.11.8.5} is satisfied,
then we have $x < 4/m_{j}^2$ from \eqref{13.19.12}, and so
\begin{equation}\label{e2.11.8.5}
x\eta < \frac{1}{4N^2m_{j}^2} < \frac{1}{4m_{j+1}^2} < a_{j}-\alpha;
\end{equation}
this follows from \eqref{e3.10.8.5}, 
\eqref{e6.2.5.5} (which implies $m_{j+1} < N m_j$) and 
\eqref{13.19.12}.  In view of the previous remarks, the proof is now complete.
\end{proof}

\section{Construction of $M_\alpha$}
\label{s10.14.6.5}
In this section, we associate to any irrational number $\alpha$,
$0<\alpha<1$, a toric $4$-manifold $M_\alpha$ of infinite topological
type.

\subsection{Notation and set-up} \label{s2.5.6.5}The two points
$\alpha$ and $\infty$ decompose $\partial \cH^2$ as
\begin{equation*}
\del\cH^2 = \del\cH^2_+\cup\del\cH^2_- \cup\{\alpha\}\cup\{\infty\},
\end{equation*}
where
\begin{equation*}
\del\cH^2_+ = \{(x,0)\in\oH{}^2: x>\alpha\},\quad
\del\cH^2_- = \{(x,0)\in\oH{}^2: x<\alpha\}.
\end{equation*}
Choose a smooth simple arc $\oZ$ in $\oH{}^2$ which joins $\alpha$ to
$\infty$ and is such that
\begin{equation*}
\oZ = Z\cup\{\alpha\}\cup\{\infty\},\quad Z\subset \cH^2.
\end{equation*}
Then $Z$ decomposes $\cH^2$ as a disjoint union
\begin{equation*}
\cH^2 = D_+\cup Z\cup D_-,
\end{equation*}
where $D_\pm$ contains $\del \cH^2_{\pm}$ in its closure. Finally, put
\begin{equation*}
\oD_{\pm} = D_{\pm}\cup \del \cH^2_{\pm}.
\end{equation*}
The reader is urged to note that $\oD_+$ does not contain $Z$ or
either of the points $\alpha$ and $\infty$.

\subsection{Notation} The boundary component $\del \cH^2_+$ is 
decomposed into the intervals $[a_j,a_{j-1}]$.  We shall refer to
these as {\em edges} and to the $a_j$ themselves as {\em corners}.

\subsection{Construction of $M_\alpha$} Recall the combinatorial
description of smooth toric $4$-manifolds of Orlik and Raymond
\cite{OrRa:at4}. Let $M$ be a compact, simply connected $4$-manifold, with a
smooth action of the 2-torus $T^2$, free on the open subset $U\subset M$.  Let
$\pi:M \to M/T^2 = P$ be the quotient map. Then $P$ is a topological polygon
with a finite number of edges $E_j$, say.  If $x$ is a point in the interior of
$P$, then $\pi^{-1}(x) = T^2$, while if $x\in \mathrm{int}\, E_j$, then
$\pi^{-1}(x)$ is a circle.  This circle is a special orbit of the $T^2$-action,
with isotropy group $S^1_{\smash{(m_j,n_j)}}$. If $x$ is a corner $E_j\cap
E_{j+1}$ of $P$, then $\pi^{-1}(x)$ is a point, fixed by the $T^2$-action.

More important for our purposes is the converse construction: to a polygon $P$,
with edges $E_j$, labelled by coprime pairs $(m_j,n_j)$, one can construct a
toric $4$-orbifold $M$ such that $M/T^2$ is $P$, with special orbits described
in the previous paragraph; for $M$ to be a smooth manifold, we require that the
pairs $(m_j,n_j)$ and $(m_{j+1}, n_{j+1})$ form a $\Z$-basis for $\Z\oplus \Z$
for each $j$. One way to understand this construction is as follows. Starting
from the closed polygon $P$, form the product $P\times T^2$. For each point $x
\in \partial P$, the labelling gives us a subtorus $T_x$ of $T^2$, which is a
circle if $x$ is in the interior of an edge and is $T^2$ itself if $x$ is a
corner. The manifold $M$ is formed by contracting $T_x$ to a point. A smooth
local model for these `contractions' is given by the toric description of
$\R^4$.  Using polar coordinates, the $T^2$-action is
\begin{equation*}
(z_1,z_2)\cdot (r_1e^{i\theta_1}, r_2e^{i\theta_2}) = 
(r_1e^{i(\theta_1+z_1)}, r_2e^{i(\theta_2+z_2)}).
\end{equation*}
The quotient space is the closed quadrant $Q = \{(r_1,r_2): r_1\geq 0,
r_2\geq 0\}$. We have maps
\begin{equation}\label{e1.20.6.5}
Q \times T^2 \stackrel{\beta}{\longrightarrow}
\R^4\stackrel{\pi_c}{\longrightarrow} Q
\end{equation}
where $\pi_c$ is the quotient map $R^4\to\R^4/T^2=Q$ and
\begin{equation*}
\beta(r_1,r_2,\theta_1,\theta_2) = (r_1e^{i\theta_1},
r_2e^{i\theta_2})
\end{equation*}
contracts $S^1_{0,1}$ along the $r_1$-axis, $S^1_{1,0}$ along the $r_2$-axis,
and the whole $T^2$ at the corner of $Q$.  Returning to the construction of
$M$, if $x\in \partial P$ is in the interior of an edge, with isotropy group
$S^1_{(m,n)}$, then we can change basis in the lattice so that $(m,n)$ is
mapped to $(0,1)$. Then a neighbourhood $N$ of $x$ in $P$ can be identified
with a neighbourhood $N'$ of $r_1=1$, say, on the $r_1$-axis in $Q$. The part
of $M$ over $N$ can then be defined to be $\pi_c^{-1}(N')$.  Similarly, if
$x$ is a corner of $P$, then we can change basis in the lattice so that the
labels of the two edges through $x$ are mapped to $(1,0)$ and $(0,1)$. Again,
we identify a small neighbourhood $N$ of $x$ with a neighbourhood $N'$ of $0$
in $Q$ and construct the part of $M$ over $N$ as $\pi^{-1}_c(N')$. This
concludes our outline of the reconstruction of a toric $4$-manifold from a
labelled polygon.

With this understood, we can describe $M_\alpha$ using the notation introduced
in \S\ref{s2.5.6.5}. That is, we form the product $\oD_+\times T^2$, and
contract $S^1_{(m_j,n_j)}$ over each point of $[a_j,a_{j-1}]$; the result is an
open manifold that we shall call $M_\alpha$.  Notice that a sequence of points
$p_j\in M_\alpha$ with $\pi(p_j)$ converging to a point on $\oZ$ will have no
convergent subsequence in $M_\alpha$. On the other hand, if $\pi(p_j)$
converges to a point on $\del \cH^2_+$, then $p_j$ {\em does} have a convergent
subsequence.

\subsection{Remark} Note that for $j\geq 1$, $S_j = \pi^{-1}[a_j,a_{j-1}]$
is a smoothly embedded $2$-sphere. Suitably oriented, we have
\begin{equation*}
S_j\cdot S_{j+1} = - 1,\quad S_j\cdot S_j = e_j.
\end{equation*}

\section{A self-dual Einstein metric on $M_\alpha$}
\label{s11.14.6.5}

For background to the material presented in this section, the reader is
referred to \cite{CaPe:emt} and \cite{CaSi:emcs}. We present the essential
formulae here: the proofs can be found in the cited papers.

Let $F(x,y)$ be a solution of the equation
\begin{equation}\label{e3.9.8.5}
  \Delta F = \tfrac{3}{4}F
\end{equation}
where $\Delta = y^2(\del_x^2 + \del_y^2)$ is the Laplacian of
$\cH^2$. Let
\begin{equation*}
f(x,y) = \sqrt{y}F(x,y)
\end{equation*}
and introduce the auxiliary quantities
\begin{equation*}
v_1 = (f_y,xf_y -y f_x),\quad v_2 = (f_x, xf_x +yf_y -f)
\end{equation*}
and
\begin{equation*}
w = f_x^2 + f_y^2 - y^{-1}ff_y,  \mbox{ so that }\ve(v_1,v_2) = yw.
\end{equation*}
(Recall from \eqref{e1.9.8.5} that $\ve$ stands for the standard
symplectic form on $\R^2$.)
Set
\begin{equation*}
g = \frac{|w|}{f^2}\left(\rd x^2 +\rd y^2 + \frac{\ve(v_1,\rd z)^2
+\ve(v_2,\rd z)^2}{w^2}\right).
\end{equation*}
This is a self-dual Einstein metric on $U\times T^2$, where $U\subset
\cH^2$ is the set on which $w\not=0$, $f\not=0$. Moreover, the sign of
the scalar curvature is opposite to the sign of $w$.

The next task is to choose the eigenfunction $F$ so as to obtain a
metric that extends smoothly over the special $T^2$-orbits. This is
very conveniently done using the following integral formulae which
show how to recover $F$ from its (renormalized) boundary value $u(x)=f(x,0)$.

\subsection{Integral formulae}
 Let $u$ be a distribution on
$\R$, viewed as the finite part of $\del \cH^2$. If $u$ does not
grow too fast at $\infty$, we define
\begin{equation}\label{e11.3.5.5}
f(x,y) = \int \frac{y^2}{2((x-x_1)^2 +y^2)^{3/2}}u(x_1)\,\rd x_1 ,
\end{equation}
which we shall also write as
\begin{equation} \label{e12.3.5.5}
f(x,y) = [k_y*u](x,y), \mbox{ where } k_y(x) = \frac{y^2}{2(x^2+y^2)^{3/2}}.
\end{equation}
Then $F(x,y) = y^{-\hf}f(x,y)$ satisfies \eqref{e3.9.8.5} and has
boundary data given by $u$ in the sense that
\begin{equation*}
f(x,y)\to u(x)\mbox{ as }y\to 0
\end{equation*}
in the sense of distributions.

If $f$ is given by \eqref{e11.3.5.5}, there is also a formula for $w$
in terms of the boundary value $u$. 
For this, introduce
\begin{equation*}
(\mu(x),\nu(x)) = (u'(x), xu'(x) -u(x)),
\end{equation*}
so that $(\mu,\nu)$ is the boundary value of $v_2$ and
$\mu(x)x-\nu(x) = u(x)$. 
If
\begin{equation} \label{e2.3.6.5}
D(x_1,x_2) = (\mu(x_1)\nu(x_2)-\mu(x_2)\nu(x_1))(x_1-x_2),
\end{equation}
then
\begin{equation}\label{e11.6.5.5}
w(x,y) =\hf y^{-2}\int\!\!\int k_y(x-x_1)k_y(x-x_2)D(x_1,x_2)\,\rd x_1 \rd x_2.
\end{equation}
This is essentially \cite[(5.16)]{CaSi:emcs}.

\subsection{Remark} At the beginning of this section we made the
assumption that $u$ should not grow too fast at $\infty$. Our
explicit choice of coordinates obscures the invariance of the
preceding formulae. It is more natural to interpret $u$ as a
distributional section of the $-\hf$-power of the density bundle on
$\del \cH^2$. Then the above integral formulae become fully invariant
under $PSL_2(\R)$, acting by isometries on $\cH^2$ and projectively on
its boundary. The distributions that we shall 
actually use will satisfy $u(\pm x) = \pm 1$ for all sufficiently
large $|x|$; it is not difficult to check that $u(x)|\rd
x|^{-\hf}$ is then smooth in a neighbourhood of $\infty$.

\subsection{Definition of $g_\alpha$}
In order to obtain a metric on $M_\alpha$, apply the previous formulae with
$u$ equal to  the odd extension of 
$\eta$ to the left of $\alpha$:
\begin{equation*}
u(x) = \eta(x) - \eta(2\alpha-x).
\end{equation*}
It is clear from \eqref{e11.3.5.5} that $f(\alpha,y)=0$ for all $y\geq
0$;  accordingly we set 
$Z = \{(\alpha,y):y>0\}$ (cf.\ \S\ref{s2.5.6.5}).

\begin{thm} The metric $g_\alpha$ is defined on $D_+\times T^2$ and
  extends smoothly to $M_\alpha$. 
\end{thm}
\begin{proof}   See \cite[Theorem 5.2.1]{CaSi:emcs}.  That result
applies to show that $w>0$ in $\cH^2$ and that $f>0$ in $D_+$, so that
$g_\alpha$ is defined on $D_+\times T^2$. We explain why the metric extends to
$M_\alpha$ since this was not done in \cite{CaSi:emcs} and we did not
explicitly check that $w$ has the needed boundary behaviour.  This is now fixed
in Proposition~\ref{p1.5.7.5} below.

For this, we need to understand the asymptotic behaviour of $f$ and $w$ as
$y\to 0$.  So pick $x_0>\alpha$. We distinguish two cases according to whether
$x_0$ is or is not one of the $a_j$. The easier case corresponds to $a_j < x_0
< a_{j-1}$ for some $j$. Then by Proposition~\ref{p1.6.6.5}, with $m_1=m_2 =
m_j$, $n_1=n_2=n_j$, we have
\begin{equation}\label{e2.5.7.5}
f(x,y) = m_jx-n_j +y^2f_1(x) +\cdots,
\end{equation}
and by part (i) of Proposition~\ref{p1.5.7.5}, 
\begin{equation}\label{e1.5.7.5}
w(x,y) = w_0(x) + w_1(x)y^2 +\cdots, \mbox{ where }w_0(x)>0,
\end{equation}
for $|x-x_0|$ and $y\geq 0$ sufficiently small.
Write $(m_j,n_j) = (m,n)$, to simplify the notation. Then computing with
these formulae,
\begin{align*}
v_1 &= y(2f_1,2xf_1 - m) + O(y^3)\\
v_2 &= (m,n) + O(y^2)\\
\ve(v_1,v_2) &= y(2uf_1 -(u')^2).
\end{align*}
Thus, for small $y$,
\begin{equation*}
g_\alpha \simeq
\frac{w_0}{u^2}\left(
\rd x^2 +\rd y^2 + \frac{
y^2(2f_1\rd z_2 + (u'-2xf_1)\rd z_1)^2
+ (m\rd z_2 -n\rd z_1)^2}{(2uf_1 - (u')^2)}\right)
\end{equation*}
In particular
\begin{equation*}
g_\alpha(m\del_{z_1}+n\del_{z_2},
m\del_{z_1}+n\del_{z_2}) \simeq y^2.
\end{equation*}
If we introduce new coordinates
\begin{align*}
\rd \theta &= m\rd z_2 - n\rd z_1 \\
\rd \psi &= m_1\rd z_2 - n_1\rd z_1,
\end{align*}
where $|mn_1-m_1n|=1$, we have
\begin{equation*}
g_\alpha = \frac{w_0}{u^2}
\left( \rd x^2 +  a(x)\rd\theta^2 + \rd y^2
+ y^2\rd \psi^2+\cdots)\right),\quad (a(x)>0)
\end{equation*}
which shows that $g_\alpha$ does extend as a smooth metric to
$\pi^{-1}(U)$, where $U$ is a small neighbourhood of the
boundary-point $x_0$. 

To complete the proof, we must consider the behaviour of $g_\alpha$ at
one of the corners $a_j$. By a suitable change of variables, we
may assume that $x_0 = a_0 = 1$, so that $u(x) =x$ for $1-\delta<x
\leq 1$, $u(x) =1$ for $1\leq x < 1+\delta$.

By Proposition~\ref{p1.6.6.5}, we have 
\begin{equation*}
f(x,y) =1 + \hf(x-1) - \hf((x-1)^2+y^2)^{1/2} + f_1(x)y^2 +\cdots
\end{equation*}
and by part (ii) of Proposition~\ref{p1.5.7.5} 
\begin{equation}\label{e4.5.7.5}
w(x,y) = \hf((x-1)^2+y^2)^{-1/2} + w_1(x) + \cdots
\end{equation}
where the omitted terms contain only even powers of $y$. Using the
change of variable
\begin{equation*}
x-1 + iy = (r_2 + ir_1)^2
\end{equation*}
we compute, to leading order,
\begin{align*}
v_1 &= -\frac{r_1r_2}{r_1^2+r_2^2}(1,1) + \cdots\\
v_2 &= \frac{1}{r_1^2+r_2^2}(r_1^2,-r_2^2) + \cdots\\
w &= \frac{1}{2}\frac{1}{r_1^2+r_2^2}+\cdots\\
\rd x^2 + \rd y^2 &=4(r_1^2+r_2^2)(\rd r_1^2 +\rd r_2^2).
\end{align*}
Substituting these into the formula for $g_\alpha$ we find
\begin{equation*}
g_\alpha = (\rd r_1^2 +r_1^2 \rd z_1^2) +
(\rd r_2^2 +r_2^2 \rd z_2^2) +\cdots
\end{equation*}
exactly as required for smooth extension to $M_\alpha$ in a
neighbourhood of the corner.
\end{proof}

\section{Boundary behaviour of $f$ and $w$}\label{s1.5.7.5}

This section is devoted to a study of the integral formulae for $f$
and $w$, \eqref{e11.3.5.5} and \eqref{e11.6.5.5}.  As in the previous
section, we assume that $u$ 
is a distribution on $\R$ such that
\begin{equation} \label{e2.6.6.5}
\mbox{$u$ is equal to a constant multiple of $\sign(x)$ outside a compact set.}
\end{equation}
Since
\begin{equation} \label{e3.6.6.5}
K_y ''(x) = k_y(x),\quad K_y(x) =\hf \sqrt{x^2+y^2},
\end{equation}
we can rewrite \eqref{e11.3.5.5} as
\begin{equation} \label{e4.6.6.5}
f(x,y) = \int K_y(x-x_1)u''(x_1)\,\rd x_1,
\end{equation}
the condition \eqref{e2.6.6.5} being used to justify the integration by parts.

We now start our study of the behaviour of \eqref{e11.3.5.5} as $y\to
0$.  For this, fix $\delta >0$, and define
\begin{equation}\label{e5.6.6.5}
B_0 = (-\delta,\delta) \times (0,\delta),\quad
B = (-\delta,\delta) \times [0,\delta).
\end{equation}
Fix also a smooth cut-off function $\beta$, equal to $1$ for $|x|\leq
2\delta$ and equal to $0$ for all $|x|\geq 3\delta$.
\begin{lem} \label{l1.6.6.5}
\textup{(i)} In \eqref{e11.3.5.5}, suppose that $u=0$ in $(-2\delta,2\delta)$.
Then $f(x,y)$ is real-analytic for $(x,y)\in B$ and has an expansion in this
domain of the form
\begin{equation}\label{e6.6.6.5}
f(x,y) = y^2f_2(x) + y^4f_4(x) + \cdots
\end{equation}
\textup{(ii)} In \eqref{e11.3.5.5}, suppose that $u$ has compact
support, and that $u$ is $C^2$ in $(-3\delta,3\delta)$. Then $f(x,y)$
is continuous for $(x,y)$ in $B$ and real-analytic for $(x,y)\in B_0$.
Moreover, if $(x,y)\in B$,
\begin{equation}\label{e7.6.6.5}
|f(x,y) - u(x)| \leq \hf |y| \int |(\beta u)''(x_1)|\,\rd x_1\mbox{ as }y\to 0.
\end{equation}
\end{lem}
\begin{proof}
For (i), note that if $(x,y)\in B$ and 
$x_1$ is in the support of $u$, we have $|x-x_1|>  \delta$ and so $k_y(x-x_1)$
can be expanded as a convergent series
\begin{equation}\label{e21.3.5.5}
\frac{y^2}{2((x-x_1)^2+y^2)^{3/2}}= 
\frac{y^2}{2|x-x_1|^3}\left(1-\frac{3}{2}\frac{y^2}{(x-x_1)^2}+\cdots \right).
\end{equation}
The result follows at once from this, using term-by-term integration.

For part (ii), note first that
\begin{equation}\label{e85.6.6.5}
f(x,y) = \int k_y(x-x_1)\beta(x_1)u(x_1)\,\rd x_1 +
\int k_y(x-x_1)(1-\beta(x_1))u(x_1)\,\rd x_1 = I_1 + I_2.
\end{equation}
By part (i), $I_2$ is real-analytic and $O(y^2)$ if $(x,y)\in B$. On the other
hand, by \eqref{e4.6.6.5} 
\begin{equation}\label{e9.6.6.5}
I_1(x,y) = \int \hf\sqrt{(x-x_1)^2+y^2}(\beta u)''(x_1)\,\rd x_1.
\end{equation}
Hence $I_1$ is real-analytic in $B_0$ and continous in $B$, with
\begin{equation}\label{e10.6.6.5}
I_1(x,0) = \int \hf |x-x_1|(\beta u)''(x_1)\,\rd x_1 = u(x).
\end{equation}
The estimate \eqref{e7.6.6.5} is obtained by noting
\begin{equation}\label{e11.6.6.5}
|I_1(x,y) - I_1(x,0)| \leq \hf
\int \bigl(\sqrt{(x-x_1)^2+y^2} - |x-x_1|\bigr)|(\beta u)''(x_1)| \,\rd x_1
\end{equation}
and using the triangle inequality,
\begin{equation*}
0\leq  \sqrt{(x-x_1)^2+y^2} - |x-x_1| \leq y.
\end{equation*}
\end{proof}

We turn to now to the boundary behaviour of \eqref{e11.3.5.5} when $u$
is piecewise linear near $0$. 

\begin{lem} Suppose in the above that
\begin{equation}\label{e12.6.6.5}
u(x) = a\,\sign(x) + b|x| + c x\mbox{ for }|x| < 3\delta.
\end{equation}
Then for $(x,y)\in B_0$,
\begin{equation}\label{e22.3.5.5}
f(x,y) = a\frac{x}{(x^2+y^2)^{1/2}} + b(x^2+y^2)^{1/2} + cx +  O(y^2)
\end{equation}
where $O(y^2)$ stands for a real-analytic function of $(x,y)$ which
goes to $0$ like $y^2$ as $y\to 0$.
\label{l1.5.7.5}\end{lem}

\begin{proof} With $\beta$ as before,
\begin{equation*}
(\beta u)''(x) = 2a\delta'(x) + 2 b \delta(x) + \beta'' u
\end{equation*}
so that, using \eqref{e4.6.6.5}, 
\begin{equation*}
f(x,y) = 2a K_y(x) + 2b K'_y(x) + \int K_y(x-x_1)\beta''(x_1)u(x_1)\,\rd x_1
= f_1(x,y) + f_2(x,y) +f_3(x,y).
\end{equation*}
Since $\beta''(x_1)u(x_1)$ vanishes for $|x|\leq 2\delta$, it follows
as in Lemma~\ref{l1.6.6.5} that $f_3(x,y)$ is real-analytic in $B$. By
part (ii) of the same lemma
 $f_3(x,0) =cx$ (for we know that
$f(x,0) = u(x)$, at least if $x\not=0$). Since
\begin{equation*}
f_1(x,y) = a\sqrt{x^2+y^2}, \quad f_2(x,y) = \frac{b x}{\sqrt{x^2+y^2}},
\end{equation*}
the proof is complete.
\end{proof}

Using these lemmas, the proof of the following result is immediate.
\begin{prop} Suppose that 
\begin{equation*}
u(x) = \begin{cases} m_1 x - n &\mbox{for }0<x< 3\delta;\\
m_2 x - n &\mbox{for }-3\delta<x< 0.
\end{cases}
\end{equation*}
Then if \eqref{e5.6.6.5} is satisfied, we have
\begin{equation*}
f(x,y) = \hf(m_1-m_2)\sqrt{x^2+y^2} + \hf(m_1+m_2)x - n + O(y^2).
\end{equation*}
\label{p1.6.6.5}
\end{prop}
\begin{proof}
Write
\begin{equation*}
u(x) = \hf(m_1-m_2)|x| + \hf(m_1+m_2)x -n
\end{equation*}
and apply Lemma~\ref{l1.5.7.5}.
\end{proof}

We now give our result concerning the boundary behaviour of
$w(x,y)$. Let the notation be as in \S\ref{s2.5.6.5}.

\begin{prop} \textup{(i)} Let $c$ be a point at which $f(x)$ is smooth. Then
$w(x,y)$ is smooth and positive in $B$.

\noindent\textup{(ii)} Let $c$ be one of the $a_j$. Then in $B_0$ we have
\begin{equation*}
w(x,y) = \frac{1}{\sqrt{(x-c)^2+y^2}} +O(1),
\end{equation*}
where $O(1)$ stands for a smooth positive function in $B$.
\label{p1.5.7.5}\end{prop}

\begin{proof}  In case (i), write  \eqref{e11.6.5.5} in the form
\begin{equation}\label{e1.27.7.5}
w(x,y) = \frac{1}{4}\int\left[ k_y(x-x_1) \int
  ((x-x_2)^2+y^2)^{-3/2}D(x_1,x_2)\,\rd x_2\right]\, \rd x_1.
\end{equation}
Because $D(x_1,x_2)=0$ if $x_1$ and $x_2$ are sufficiently close to
$c$, the function
\begin{equation}\label{e2.27.7.5}
x\mapsto \int |x-x_2|^{-3}D(x,x_2)\,\rd x_2
\end{equation}
is smooth for $x$ near $c$, and because $k_y(x) \to \delta(x)$ as
$y\to 0$, we obtain
\begin{equation*}
w(x,0) = \frac{1}{4}\int |x-x_2|^{-3}D(x,x_2)\,\rd x_2
\mbox{ for $x$ sufficiently   close to $c$.}
\end{equation*}
In particular, $w(x,0)>0$ by the non-negativity of $D(x_1,x_2)$. With
a little more work it can be shown that $w(x,y)$ is smooth in $B$,
following the proof of Lemma~\ref{l1.6.6.5}.

In case (ii), we use a cut-off function $\beta$, equal to $1$ near $c$
and with small support, to split the integral as
\begin{equation*}
w(x,y) = w_1(x,y) + w_2(x,y),
\end{equation*}
where
\begin{equation*}
w_1(x,y) = \frac{1}{4}\int\left[ k_y(x-x_1) \int
  ((x-x_2)^2+y^2)^{-3/2}(1-\beta(x_2))D(x_1,x_2)\,\rd x_2\right]\,
\rd x_1,
\end{equation*}
and
\begin{equation*}
w_2(x,y) = \frac{1}{4}\int\left[ k_y(x-x_1) \int
  ((x-x_2)^2+y^2)^{-3/2}\beta(x_2)D(x_1,x_2)\,\rd x_2\right]\,
\rd x_1.
\end{equation*}
Now $w_1(x,0)$ is smooth and positive near $c$ by the same argument as
before, and if we write
\begin{equation*}
w_2(x,y) = w_3(x,y) + w_4(x,y)
\end{equation*}
where
\begin{equation*}
w_3(x,y) = \frac{1}{4}\int\left[ k_y(x-x_1) (1-\beta(x_1))\int
  ((x-x_2)^2+y^2)^{-3/2}\beta(x_2)D(x_1,x_2)\,\rd x_2\right]\,
\rd x_1,
\end{equation*}
and
\begin{equation*}
w_4(x,y) = \frac{1}{4}\int\left[ k_y(x-x_1) \beta(x_1)\int
  ((x-x_2)^2+y^2)^{-3/2}\beta(x_2)D(x_1,x_2)\,\rd x_2\right]\,
\rd x_1,
\end{equation*}
then $w_3(x,0)$ is also smooth and positive near $c$ by symmetry. Thus
it remains only to consider $w_4(x,y)$.

In order to simplify the notation, use a translation to set $c=0$ and
assume that  $u$ is as in
Proposition~\ref{p1.6.6.5} near $x=0$. Then
\begin{equation*}
D(x_1,x_2) = 
-\hf x_2\sign(x_1) -\hf x_1\sign(x_2) + \hf|x_1| + \hf|x_2|
\mbox{ for }(x_1,x_2) \in (-\delta,\delta)\times (-\delta,\delta),
\end{equation*}
assuming, as we may, that
\begin{equation*}
n(m_1-m_2) =1.
\end{equation*}

Now Lemma~\ref{l1.5.7.5} can be applied, giving
\begin{equation*}
\int k_y(x-x_1)\beta(x_1)D(x_1,x_2)\,\rd x_1
=\hf
- \frac{xx_2}{\sqrt{x^2+y^2}} + \sqrt{x^2+y^2}
- x_1 \sign(x_2) + |x_2| + O(y^2)
\end{equation*}
and then
\begin{align*}
\int\!\!\int k_y(x-x_1)\beta(x_1)\beta(x_2)D(x_1,x_2)\,\rd x_1\rd x_2 
&= 
- \frac{x^2}{\sqrt{x^2+y^2}} + \sqrt{x^2+y^2} + O(y^2) \\
&= \frac{y^2}{\sqrt{x^2+y^2}} + O(y^2).
\end{align*}
Dividing this by $y^2$ and combining with our conclusions about the
other $w_j$ now gives the stated result.
\end{proof}

\section{Completeness of $g_\alpha$}
We now have a SDE metric on $M_\alpha$. Our final task is to show that
this metric is complete. 
\label{s1.9.8.5}
\begin{thm}
Let
$\Gamma:[0,1)\to M_\alpha$ be a smooth curve with finite length,
\begin{equation}\label{e1.1.12.4}
L(\gamma) =  \lim_{t\to 1} d(\Gamma(0),\Gamma(t)) < \infty.
\end{equation}
Then there exists a point $p\in M_\alpha$ such that
\begin{equation}\label{e1.6.5.5}
\lim_{t\to 1} \Gamma(t) = p.
\end{equation}
\end{thm}

\begin{proof}
Let $\gamma = \pi\circ \Gamma$ be the projection of the curve to
$\oD_+$. Then from the form of the metric $g_\alpha$, it is clear that this
too has finite length with respect to the base metric
\begin{equation}\label{3.22.12}
h = \frac{w} {f^2}(\rd x^2 + \rd y^2).
\end{equation}

From the estimates proved below, $h$ is uniformly bounded below by a multiple
of the euclidean metric on $\R^2$, so that
\begin{equation*}
\gamma(1) : =  \lim_{t\to 1}\gamma(t)
\end{equation*}
exists and lies in $\oD_+\cup Z\cup\{\alpha\}$.  There are now three
possibilities to consider. First, suppose $\gamma(1)\in \oD_+$. Then
$\pi^{-1}(\gamma(1))$ is contained in the interior of $M_\alpha$, and
so \eqref{e1.6.5.5} must be satisfied. Next, suppose, if possible,
that $\gamma(1)\in Z$. Since $f$ vanishes on $Z$, $h$ has a
double-pole along $Z$, and it follows that $\gamma$ must have infinite
length, a contradiction.  The remaining possibility is that
$\gamma(1) = (\alpha,0)$.

From the estimates proved below, we have
\begin{equation}\label{e3.6.5.5}
x>\alpha, y>0\mbox{ implies that }\frac{w}{f^2}>C(x+y-\alpha)^{-2}.
\end{equation}
If we set $\xi =x -\alpha - y$, $\eta = x -\alpha +y$, then we get
\begin{equation*}
h \geq C\frac{\rd \xi^2 + \rd \eta^2}{2\eta^2}\mbox{ for }|\xi| \leq \eta
\end{equation*}
(that is, LHS minus RHS is positive-definite).  Now our curve $\gamma$
is contained in $\{|\xi| < \eta\}$ and $\lim_{t\to 1}\gamma(t) =
(0,0)$. Hence $l(\gamma)= \infty$ and this contradiction completes the
proof, modulo the estimates established in the next section.
\end{proof}

\subsection{Estimates}

In order to simplify the notation in this section, we shift variables
so that $\alpha$ is translated to the origin, and we aim to understand the
behaviour of $f$ near $(0,0)$.

First of all, we have
\begin{prop}\label{p1.5.3.4}
Let $u(x)$ and $f(x,y)$ be as throughout.  Then there exist $\ve>0$ and $C>0$
so that
\begin{equation}\label{e4.6.5.5}
f(x,y)\leq  C\sqrt{x+y} \mbox{ if } (x,y) \in (0,\ve)\times (0,\ve).
\end{equation}
\end{prop}
\begin{proof}
The result will be established by showing that
\begin{equation}\label{e1.22.8.5}
f(x,\theta x) \leq  C\sqrt{x}\mbox{ if }0<\theta \leq 1
\end{equation}
and
\begin{equation}\label{e2.22.8.5}
f(\theta_1 y, y) \leq  C\sqrt{y}\mbox{ if }0<\theta_1 \leq 1.
\end{equation}
Since $u$ is odd,
\begin{equation}\label{e1.3.6.5}\begin{split}
f(x,y) &=\int_{0}^\infty \{k_y(x-x_1) - k_y(x+x_1)\}u(x_1)\,\rd x_1\\
&\leq \Omega\int_{0}^\infty \{k_y(x-x_1) - k_y(x+x_1)\}\sqrt{x_1}\,\rd x_1\\
&\leq \Omega\int_{0}^\infty k_y(x-x_1)\sqrt{x_1}\,\rd x_1,
\end{split}
\end{equation}
where we have used the upper bound of Theorem~\ref{t1.3.5.5} and the
positivity of $k_y(x)$. 

To prove \eqref{e1.22.8.5}, introduce a natural rescaling of variables
in \eqref{e1.3.6.5}, 
\begin{equation}\label{e7.6.5.5}
\xi = x_1/x, \theta = y/x
\end{equation}
so that \eqref{e1.3.6.5} becomes
\begin{equation} \label{1.27.5}
f(x,x\theta) \leq \Omega I = \Omega\sqrt{x}\int_0^\infty
k_\theta(\xi-1)\sqrt{\xi}\,\rd \xi. 
\end{equation}
Note that this yields \eqref{e1.22.8.5} for each fixed positive
$\theta$, but the uniformity as $\theta\to 0$ needs a little further work.
For this   we split $I$, as in \S\ref{s1.5.7.5}, using a bump-function
$\beta$ identically 
equal to $1$ in a neighbourhood of $\xi=1$. We have 
\begin{equation*}
I = I_1 + I_2
\end{equation*}
where
\begin{equation*}
I_1 =\int_0^\infty k_\theta(\xi-1)(1-\beta(\xi))\sqrt{\xi}\,\rd \xi,\quad
I_2 =\int_0^\infty k_\theta(\xi-1)\beta(\xi)\sqrt{\xi}\,\rd \xi.
\end{equation*}
Then since $\lim_{\theta\to 0}k_\theta(\xi-1) =\delta(\xi-1)$, we have
\begin{equation*}
I_2\to  Cx^{1/2}\beta(1)\sqrt{1}= C\sqrt{x} \mbox{ as }\theta \to 0.
\end{equation*}
By continuity, $I_2$ is uniformly bounded by $C\sqrt{x}$ for all
$0\leq \theta \leq 1$. The estimate \eqref{e1.22.8.5} now follows by
noting that  $I_1\geq 0$ and $I_1\to 0$ as $\theta\to
0$.  In fact, it is easily
shown that $I_1 \leq C' \theta^2\sqrt{x}$ for
$0\leq\theta \leq 1$.

The complementary estimate \eqref{e2.22.8.5} is somewhat simpler:
return to \eqref{e1.3.6.5} and make the change of
variables $x = y\theta_1$, $x_1 = y\xi_1$, so
\begin{equation}\label{2.27.5}
f(y\theta_1,y ) \leq \Omega\sqrt{y}
\int_0^\infty k_1(\xi_1-\theta_1)\sqrt{\xi_1}\,\rd \xi_1.
\end{equation}
Clearly the integral is uniformly bounded for $\theta_1\in [0,1]$ 
and this completes the proof.
\end{proof}

Next we need a lower bound on $w(x,y)$. This is given by
\begin{prop}\label{p2.5.3.4}
Suppose that in the continued fraction expansion of $\alpha$, $3\leq e_j \leq
N$, for some $N$.  Then there exist $\ve>0$ and $C>0$ so that
\begin{equation}\label{e3.5.3.4}
w(x,y)\geq  C(x+y)^{-1}\mbox{ if } (x,y) \in (0,\ve)\times (0,\ve).
\end{equation}
\end{prop}
\begin{proof} The argument is closely analogous to the proof of the previous
proposition. In particular, \eqref{e3.5.3.4} will be established by proving the
separate inequalities
\begin{equation}
w(x,\theta x) \geq  Cx^{-1}\mbox{ if }0<\theta \leq 1
\end{equation}
and
\begin{equation}\label{e5.9.8.5}
w(\theta_1 y, y) \geq  Cy^{-1}\mbox{ if }0<\theta_1 \leq 1.
\end{equation}

We use~\eqref{e11.6.5.5}. Because $D\geq 0$ and
$k_y(x)>0$, it is enough to prove the lower bound for
\begin{equation}\label{e3.3.6.5}
I(x,y) := \hf y^{-2}\int_{x_1=0}^\infty \int_{x_2=0}^{x_1}k_y(x-x_1)k_y(x-x_2)
D(x_1,x_2)\,\rd x_1\rd x_2
\end{equation}
where the integral has been restricted to the intersection of the
positive quadrant with the region $x_2\leq x_1$.   In this integral,
make the change of variables
\begin{equation*}
x_1 = x\xi, x_2 = x\xi\eta, y = \theta x.
\end{equation*}
Then
\begin{equation*}
I(x,\theta x) =
\frac{1}{2x^2\theta^2}\int_{\xi=0}^\infty\int_{\eta=0}^1
k_\theta(\xi-1)k_\theta(\xi\eta-1)D(x\xi,x\xi\eta)\xi\rd\xi\rd \eta
\end{equation*}
\begin{equation}\label{e5.3.6.5}
\hspace{1in} \geq
\frac{1}{2x\theta^2}\int_{\xi=0}^{a/x}\int_{\eta=0}^b
k_\theta(\xi-1)k_\theta(\xi\eta-1)\xi^2(1-\eta)\rd\xi\rd \eta,
\end{equation}
where $a = 1 -\alpha$, $b = 1/16N^2$ as in Lemma~\ref{l1.11.8.5}.  For each
$\theta$ this gives the required $O(1/x)$ lower bound. To see this is uniform
as $\theta \to 0$, rewrite \eqref{e5.3.6.5} as follows and take the limit:
\begin{equation}\label{e33.3.6.5}
\frac{1}{4x} \int_0^{a/x} k_\theta(\xi-1)\left[\xi^2
\int_{\eta=0}^b(1-\eta)((1-\xi\eta)^2 +\theta^2)^{-3/2}\rd\eta\right] \rd \xi
\longrightarrow
\frac{1}{4x}
\int_0^b(1-\eta)^{-2}\rd\eta\mbox{ as }\theta \to 0.
\end{equation}
By continuity, the required uniform lower bound follows.

In order to obtain the other bound \eqref{e5.9.8.5}, return to \eqref{e3.3.6.5}
and make the substitutions
\begin{equation*}
x = \theta_1 y, x_1 = y\xi, x_2 = y\xi\eta
\end{equation*}
to give
\begin{equation}\label{e5.11.8.5}
w(y\theta_1, y) \geq
\frac{1}{2y}\int_{\xi=0}^{a/y}\int_0^b
k_1(\xi-\theta_1)k_1(\xi\eta-\theta_1)\xi^2(1-\eta)\,\rd\xi\rd\eta.
\end{equation}
This gives \eqref{e5.9.8.5} at once.
\end{proof}

\subsection{Remark}  As in \cite[\S5]{CaSi:emcs} the distribution
$u$ can be perturbed to the left of $\alpha$ to yield an infinite-dimensional
family of SDE metrics on $M_\alpha$. These perturbations need to preserve the
monotonicity and convexity properties enjoyed by $u$, as well as the boundary
condition $u=-1$ for $x\ll 0$. Provided these perturbations are sufficiently
small and supported away from $\alpha$, the resulting metrics will be complete;
the details are left to the interested reader.

%
%
%
\newcommand{\bauth}[1]{\mbox{#1}} \newcommand{\bart}[1]{\textit{#1}}
\newcommand{\bjourn}[4]{#1\ifx{}{#2}\else{ \textbf{#2}}\fi{ (#4)}}
\newcommand{\bbook}[1]{\textsl{#1}}
\newcommand{\bseries}[2]{#1\ifx{}{#2}\else{ \textbf{#2}}\fi}
\newcommand{\bpp}[1]{#1} \newcommand{\bdate}[1]{ (#1)} \def\band/{and}
\newif\ifbibtex
\ifbibtex

\bibliographystyle{genbib}
\bibliography{papers}

\else

\fi
\end{document}
\bibitem{LeBr:hcc}
\bauth{C.~R. LeBrun}, \bart{{$\mathcal H$}-space with a cosmological constant},
  \bjourn{Proc. Roy. Soc. London}{A 380}{}{1982} \bpp{171--185}.

\bibitem{LeBr:pac}
\bauth{C.~R. LeBrun}, \bart{Counterexamples to the generalized positive action
  conjecture}, \bjourn{Comm. Math. Phys.}{118}{}{1988} \bpp{591--596}.

\bibitem{LeBr:cp2}
\bauth{C.~R. LeBrun}, \bart{Explicit self-dual metrics on
  {$\CP2\connect\cdots\connect\CP2$}}, \bjourn{J.~Diff. Geom.}{34}{}{1991}
  \bpp{223--253}.

\bibitem{LeBr:cqk}
\bauth{C.~R. LeBrun}, \bart{On complete quaternionic-{K\"a}hler manifolds},
  \bjourn{Duke Math. J.}{63}{}{1991} \bpp{723--743}.

\bibitem{RM:pc}
\bauth{R.~R. Mazzeo}, \bart{Private communication.}

\bibitem{Ped:emm}
\bauth{H.~Pedersen}, \bart{{E}instein metrics, spinning top motions and
  monopoles}, \bjourn{Math. Ann.}{274}{}{1986} \bpp{35--39}.

\bibitem{Rol:reh}
\bauth{Y.~Rollin}, \bart{Rigidit\_1e d_1{E}instein du plan hyperbolique
  complexe}, \bjourn{C. R. Math. Acad. Sci. Paris}{334}{}{2002} \bpp{671--676}.

\bibitem{Tod:p6}
\bauth{K.~P. Tod}, \bart{Self-dual {E}instein metrics from the {P}ainlev\_1e
  {VI} equation}, \bjourn{Phys. Lett.}{A 190}{}{1994} \bpp{221--224}.

\end{thebibliography}

\fi

\end{thebibliography}
\end{document}